\documentclass[12pt]{article}

\usepackage{amsmath,amssymb,amsthm,amscd,a4wide,upref}

\usepackage[enableskew]{youngtab}
\usepackage{ytableau}

\usepackage{hyperref}
\hypersetup{colorlinks,citecolor=blue,filecolor=black,linkcolor=blue,urlcolor=blue}

\usepackage{lmodern}     
\usepackage[T1]{fontenc}

\begin{document}


\newcommand{\non}{\nonumber}
\newcommand{\scl}{\scriptstyle}
\newcommand{\sclnearrow}{{\scl\nearrow}\ts}
\newcommand{\scloplus}{{\scl\bigoplus}}
\newcommand{\wt}{\widetilde}
\newcommand{\wh}{\widehat}
\newcommand{\ot}{\otimes}
\newcommand{\fand}{\quad\text{and}\quad}
\newcommand{\Fand}{\qquad\text{and}\qquad}
\newcommand{\ts}{\,}
\newcommand{\tss}{\hspace{1pt}}
\newcommand{\lan}{\langle\ts}
\newcommand{\ran}{\ts\rangle}
\newcommand{\vl}{\tss|\tss}
\newcommand{\qin}{q^{-1}}
\newcommand{\tpr}{t^{\tss\prime}}
\newcommand{\spr}{s^{\tss\prime}}
\newcommand{\di}{\partial}
\newcommand{\hra}{\hookrightarrow}
\newcommand{\antiddots}
    {\underset{\displaystyle\cdot\quad\ }
    {\overset{\displaystyle\quad\ \cdot}{\cdot}}}
\newcommand{\dddots}
    {\underset{\displaystyle\quad\ \cdot}
    {\overset{\displaystyle\cdot\quad\ }{\cdot}}}
\newcommand{\atopn}[2]{\genfrac{}{}{0pt}{}{#1}{#2}}

\newcommand{\su}{s^{}}
\newcommand{\vac}{\mathbf{1}}
\newcommand{\vacf}{\tss|0\rangle}
\newcommand{\BL}{ {\overline L}}
\newcommand{\BE}{ {\overline E}}
\newcommand{\BP}{ {\overline P}}
\newcommand{\BT}{ {\overline T}}
\newcommand{\oPhi}{ {\overline \Phi}}
\newcommand{\ovp}{ {\overline \vp}}
\newcommand{\opsi}{ {\overline \psi}}
\newcommand{\ol}{\overline}
\newcommand{\pr}{^{\tss\prime}}
\newcommand{\ba}{\bar{a}}
\newcommand{\bb}{\bar{b}}
\newcommand{\eb}{\bar{e}}
\newcommand{\bi}{\bar{\imath}}
\newcommand{\bj}{\bar{\jmath}}
\newcommand{\bk}{\bar{k}}
\newcommand{\bl}{\bar{l}}
\newcommand{\hb}{\mathbf{h}}
\newcommand{\gb}{\mathbf{g}}
\newcommand{\For}{\qquad\text{or}\qquad}
\newcommand{\OR}{\qquad\text{or}\qquad}
\newcommand{\emp}{\mbox{}}


\newcommand{\U}{{\rm U}}
\newcommand{\Z}{{\rm Z}}
\newcommand{\ZY}{{\rm ZY}}
\newcommand{\Ar}{{\rm A}}
\newcommand{\Br}{{\rm B}}
\newcommand{\Cr}{{\rm C}}
\newcommand{\Fr}{{\rm F}}
\newcommand{\Mr}{{\rm M}}
\newcommand{\Sr}{{\rm S}}
\newcommand{\Prm}{{\rm P}}
\newcommand{\Lr}{{\rm L}}
\newcommand{\Ir}{{\rm I}}
\newcommand{\Jr}{{\rm J}}
\newcommand{\Qr}{{\rm Q}}
\newcommand{\Rr}{{\rm R}}
\newcommand{\X}{{\rm X}}
\newcommand{\Y}{{\rm Y}}
\newcommand{\DY}{ {\rm DY}}
\newcommand{\Or}{{\rm O}}
\newcommand{\SO}{{\rm SO}}
\newcommand{\GL}{{\rm GL}}
\newcommand{\Spr}{{\rm Sp}}
\newcommand{\Zr}{{\rm Z}}
\newcommand{\ev}{{\rm ev}}
\newcommand{\Pf}{{\rm Pf}}
\newcommand{\Ann}{{\rm{Ann}\ts}}
\newcommand{\Norm}{{\rm Norm\tss}}
\newcommand{\Ad}{{\rm Ad}}
\newcommand{\SY}{{\rm SY}}
\newcommand{\Pff}{{\rm Pf}\tss}
\newcommand{\Hf}{{\rm Hf}\tss}
\newcommand{\trts}{{\rm tr}\ts}
\newcommand{\otr}{{\rm otr}}
\newcommand{\row}{{\rm row}}
\newcommand{\End}{{\rm{End}\ts}}
\newcommand{\Mat}{{\rm{Mat}}}
\newcommand{\Hom}{{\rm{Hom}}}
\newcommand{\id}{{\rm id}}
\newcommand{\middd}{{\rm mid}}
\newcommand{\ch}{{\rm{ch}\ts}}
\newcommand{\ind}{{\rm{ind}\ts}}
\newcommand{\Normts}{{\rm{Norm}\ts}}
\newcommand{\mult}{{\rm{mult}}}
\newcommand{\per}{{\rm per}\ts}
\newcommand{\sgn}{{\rm sgn}\ts}
\newcommand{\sign}{{\rm sign}\ts}
\newcommand{\qdet}{{\rm qdet}\ts}
\newcommand{\sdet}{{\rm sdet}\ts}
\newcommand{\Ber}{{\rm Ber}\ts}
\newcommand{\inv}{{\rm inv}\ts}
\newcommand{\inva}{{\rm inv}}
\newcommand{\grts}{{\rm gr}\ts}
\newcommand{\grpr}{{\rm gr}^{\tss\prime}\ts}
\newcommand{\degpr}{{\rm deg}^{\tss\prime}\tss}
\newcommand{\Cond}{ {\rm Cond}\tss}
\newcommand{\Fun}{{\rm{Fun}\ts}}
\newcommand{\Rep}{{\rm{Rep}\ts}}
\newcommand{\sh}{{\rm{sh}}}
\newcommand{\weight}{{\rm{wt}\ts}}
\newcommand{\chara}{{\rm{char}\ts}}
\newcommand{\diag}{ {\rm diag}}
\newcommand{\Bos}{ {\rm Bos}}
\newcommand{\Ferm}{ {\rm Ferm}}
\newcommand{\cdet}{ {\rm cdet}}
\newcommand{\rdet}{ {\rm rdet}}
\newcommand{\imm}{ {\rm imm}}
\newcommand{\ad}{ {\rm ad}}
\newcommand{\tr}{ {\rm tr}}
\newcommand{\gr}{ {\rm gr}\tss}
\newcommand{\str}{ {\rm str}}
\newcommand{\loc}{{\rm loc}}
\newcommand{\Gr}{{\rm G}}

\newcommand{\twobar}{{\bar 2}}
\newcommand{\threebar}{{\bar 3}}


\newcommand{\AAb}{\mathbb{A}\tss}
\newcommand{\CC}{\mathbb{C}}
\newcommand{\KK}{\mathbb{K}\tss}
\newcommand{\QQ}{\mathbb{Q}\tss}
\newcommand{\SSb}{\mathbb{S}\tss}
\newcommand{\TT}{\mathbb{T}\tss}
\newcommand{\ZZ}{\mathbb{Z}\tss}
\newcommand{\Sbb}{\mathbb{S}}
\newcommand{\ZZb}{\mathbb{Z}}


\newcommand{\Ac}{{\mathcal A}}
\newcommand{\Bc}{{\mathcal B}}
\newcommand{\Cc}{{\mathcal C}}
\newcommand{\Cl}{{\mathcal Cl}}
\newcommand{\Dc}{{\mathcal D}}
\newcommand{\Ec}{{\mathcal E}}
\newcommand{\Fc}{{\mathcal F}}
\newcommand{\Jc}{{\mathcal J}}
\newcommand{\Gc}{{\mathcal G}}
\newcommand{\Hc}{{\mathcal H}}
\newcommand{\Lc}{{\mathcal L}}
\newcommand{\Nc}{{\mathcal N}}
\newcommand{\Xc}{{\mathcal X}}
\newcommand{\Yc}{{\mathcal Y}}
\newcommand{\Oc}{{\mathcal O}}
\newcommand{\Pc}{{\mathcal P}}
\newcommand{\Qc}{{\mathcal Q}}
\newcommand{\Rc}{{\mathcal R}}
\newcommand{\Sc}{{\mathcal S}}
\newcommand{\Tc}{{\mathcal T}}
\newcommand{\Uc}{{\mathcal U}}
\newcommand{\Vc}{{\mathcal V}}
\newcommand{\Wc}{{\mathcal W}}
\newcommand{\Zc}{{\mathcal Z}}
\newcommand{\HC}{{\mathcal {HC}}}


\newcommand{\asf}{\mathsf a}
\newcommand{\bsf}{\mathsf b}
\newcommand{\csf}{\mathsf c}
\newcommand{\nsf}{\mathsf n}


\newcommand{\Sym}{\mathfrak S}
\newcommand{\h}{\mathfrak h}
\newcommand{\q}{\mathfrak q}
\newcommand{\n}{\mathfrak n}
\newcommand{\m}{\mathfrak m}
\newcommand{\p}{\mathfrak p}
\newcommand{\gl}{\mathfrak{gl}}
\newcommand{\oa}{\mathfrak{o}}
\newcommand{\spa}{\mathfrak{sp}}
\newcommand{\osp}{\mathfrak{osp}}
\newcommand{\g}{\mathfrak{g}}
\newcommand{\kgot}{\mathfrak{k}}
\newcommand{\agot}{\mathfrak{a}}
\newcommand{\bgot}{\mathfrak{b}}
\newcommand{\sll}{\mathfrak{sl}}
\newcommand{\f}{\mathfrak{f}}
\newcommand{\z}{\mathfrak{z}}
\newcommand{\Zgot}{\mathfrak{Z}}


\newcommand{\al}{\alpha}
\newcommand{\be}{\beta}
\newcommand{\ga}{\gamma}
\newcommand{\de}{\delta}
\newcommand{\De}{\Delta}
\newcommand{\Ga}{\Gamma}
\newcommand{\ep}{\epsilon}
\newcommand{\ee}{\epsilon^{}}
\newcommand{\ve}{\varepsilon}
\newcommand{\ls}{\ts\lambda\ts}
\newcommand{\vk}{\varkappa}
\newcommand{\vs}{\varsigma}
\newcommand{\vt}{\vartheta}
\newcommand{\ka}{\kappa}
\newcommand{\vp}{\varphi}
\newcommand{\la}{\lambda}
\newcommand{\La}{\Lambda}
\newcommand{\si}{\sigma}
\newcommand{\ze}{\zeta}
\newcommand{\om}{\omega}
\newcommand{\Om}{\Omega}
\newcommand{\up}{\upsilon}


\newtheorem{thm}{Theorem}[section]
\newtheorem{lemma}[thm]{Lemma}
\newtheorem{prop}[thm]{Proposition}
\newtheorem{cor}[thm]{Corollary}
\newtheorem{conj}[thm]{Conjecture}

\theoremstyle{definition}
\newtheorem{definition}[thm]{Definition}
\newtheorem{example}[thm]{Example}

\theoremstyle{remark}
\newtheorem{remark}[thm]{Remark}

\newcommand{\bth}{\begin{thm}}
\renewcommand{\eth}{\end{thm}}
\newcommand{\bpr}{\begin{prop}}
\newcommand{\epr}{\end{prop}}
\newcommand{\ble}{\begin{lemma}}
\newcommand{\ele}{\end{lemma}}
\newcommand{\bco}{\begin{cor}}
\newcommand{\eco}{\end{cor}}
\newcommand{\bex}{\begin{example}}
\newcommand{\eex}{\end{example}}
\newcommand{\bde}{\begin{definition}}
\newcommand{\ede}{\end{definition}}
\newcommand{\bre}{\begin{remark}}
\newcommand{\ere}{\end{remark}}
\newcommand{\bcj}{\begin{conj}}
\newcommand{\ecj}{\end{conj}}

\renewcommand{\theequation}{\arabic{section}.\arabic{equation}}

\numberwithin{equation}{section}


\newcommand{\bpf}{\begin{proof}}
\newcommand{\epf}{\end{proof}}


\def\beql#1{\begin{equation}\label{#1}}

\newcommand{\bal}{\begin{aligned}}
\newcommand{\eal}{\end{aligned}}
\newcommand{\beq}{\begin{equation}}
\newcommand{\eeq}{\end{equation}}
\newcommand{\ben}{\begin{equation*}}
\newcommand{\een}{\end{equation*}}

\title{\Large\bf Fusion procedure for the Hecke--Clifford superalgebra}

\author{{Naihuan Jing,\quad Iryna Kashuba,\quad Ming Liu,}\\[0.4em]
{Alexander Molev\quad
and\quad Jian Zhang}}

\date{} 
\maketitle


\begin{abstract}
We develop a fusion procedure for the Hecke--Clifford superalgebra to show
that its idempotents can be obtained by evaluating universal
rational function. We also give a modified fusion formula for the Hecke algebra.



\end{abstract}

\section{Introduction}
\label{sec:int}

Primitive idempotents for the symmetric group $\Sym_n$ can be naturally parameterized
by the standard Young tableaux $\Uc$ with $n$ boxes. Explicit formulas for the
corresponding idempotents $e^{}_\Uc$ were given independently by Jucys~\cite{j:fy}
and Murphy~\cite{m:nc} in terms of the sums of transpositions $(1,a)+\dots+(a-1,a)$
now commonly called the {\em Jucys--Murphy elements}.

By an earlier observation of Jucys~\cite{j:yo}, the primitive idempotents
$e^{}_\Uc$ can be obtained by taking certain
limit values of the rational function
\beql{phiuuusym}
\Phi(u_1,\dots,u_n)=
\prod_{1\leqslant a<b\leqslant n}\Big(1-\frac{(a,b)}{u_a-u_b}\Big),
\eeq
where the product is taken in
the lexicographical order on the set of pairs $(a,b)$. This way of expressing
the idempotents $e^{}_\Uc$ re-emerged later as a {\em fusion procedure}
providing new solutions of the {\em Yang--Baxter equation} from old ones;
see
Cherednik~\cite{c:sb}, \cite{c:ni} and Nazarov~\cite{n:yc}.
Namely,
the factors in \eqref{phiuuusym} evaluated in the context of the Schur--Weyl duality
coincide with the {\em Yang $R$-matrices}. A Hecke algebra counterpart of the fusion
procedure was explained in detail by Nazarov~\cite{n:mh}. Its simple version
was then given in \cite{imo:ih} by extending the approach of \cite{m:fp} for
the symmetric group.
We refer the reader e.g. to \cite[Ch.~6,7]{m:yc} for applications of the fusion formula
to representations of the general linear Lie algebra and associated Yangian
along with a brief review of the literature on the extension of the fusion procedure to other classes
of algebras.

More recently, primitive idempotents for the {\em Sergeev superalgebra} $\Sc_n$ \cite{s:ta}
were constructed independently in
\cite{kms:jm}
and \cite{ls:sb}.
The role of the standard tableaux is now played by
the {\em standard shifted barred tableaux} $\Uc$. A version of the
fusion procedure expressing the corresponding idempotents $e^{}_{\Uc}$
as limit values of a rational function analogous to \eqref{phiuuusym}
was given in \cite{kms:jm} by extending the original work of Nazarov~\cite{n:ci}.

In this paper, we give a `quantum version' of the results of \cite{kms:jm} by
proving fusion formulas for the {\em Hecke--Clifford superalgebra} $\HC_n$.
This superalgebra is a $q$-analogue of $\Sc_n$,
and it was originally introduced by Olshanski~\cite{o:qu}.
Its representation theory was developed by Jones and Nazarov~\cite{jn:as},
Brundan and Kleshchev~\cite{bk:hc} and, more recently, by Li and Shi~\cite{ls:sb}.

We introduce two counterparts $\Phi(u_1,\dots,u_n)$ and $\oPhi\tss(u_1,\dots,u_n)$ of the rational function \eqref{phiuuusym} for the superalgebra $\HC_n$. Our main result provides fusion formulas
for the idempotents $e^{}_{\Uc}$ parameterized by the standard shifted barred tableaux
with $n$ boxes. They are obtained by consecutive evaluations of the functions
$\Phi(u_1,\dots,u_n)$ and $\oPhi\tss(u_1,\dots,u_n)$ at the $q$-contents of $\Uc$
(Theorems~\ref{thm:fpheckecli} and \ref{thm:fpheckeclianti}).

We believe that the fusion procedure for the superalgebra $\HC_n$ will be useful
for the construction of the associated {\em quantum immanants} and proving {\em higher Capelli identities}
in the same way as the fusion procedure for the symmetric group was essential
in the original constructions by Okounkov~\cite{o:qi} for the general linear Lie algebra;
cf. the Sergeev superalgebra version \cite{km:qi}.

In our proofs, we modify the arguments of \cite{imo:ih} and adapt them to the super-context.
We begin with a proof of the modified fusion formula for the Hecke algebra in Sec.~\ref{sec:mf}
(see Theorem~\ref{thm:fphecke}). We show that it is equivalent to the results of \cite{imo:ih}
in Proposition~\ref{prop:equivph}.
We then use a similar approach for the Hecke--Clifford superalgebra in Sec.~\ref{sec:fp}.

\section{Modified fusion formula for the Hecke algebra}
\label{sec:mf}

Fix a nonzero complex number $q$ which is not a root of unity.
The
{\em Hecke algebra} $\Hc_n$ over $\CC$ is generated by elements $T_1,\dots,T_{n-1}$
subject to the defining relations
\beql{defhecke}
\begin{aligned}
T_a \, T_{a+1} \, T_a &=
T_{a+1} \, T_a \,  T_{a+1},\\[0.4em]
T_{a} \,  T_{b} &= T_{b} \,  T_{a}\quad
{\rm for} \quad |a-b| > 1,
\\[0.4em]
T^2_a &= 1 + \ve\ts  T_a,
\end{aligned}
\eeq
where we set $\ve = q-\qin$.
For $a=1, \dots, n-1$
we let $s_{a}=(a, a+1)$ be the adjacent transpositions in the symmetric group $\Sym_{n}$.
Choose a reduced decomposition
$s=s_{a_{1}} \dots s_{a_{l}}$ of any element
$s \in \Sym_{n}$ and set $T_{s}=T_{a_{1}} \dots T_{a_{l}}$.
These elements do not depend on the chosed
reduced decompositions and form a basis of $\Hc_{n}$.
For any $a=1,\dots,n$ we let $\om_a$ denote the
unique longest element
of the group $\Sym_a$ which is regarded as
a natural subgroup of $\Sym_n$.

The irreducible representations of $\Hc_{n}$ over $\CC$ are parameterized by partitions $\la$
of $n$. Such a
partition
$\la=(\la_{1}, \dots, \la_{l})$ will be identified with the corresponding Young diagram
which is a left-justified array of rows of unit boxes, where the top row
contains $\la_1$ boxes, the next row contains $\la_2$ boxes, etc.
We will need the conjugate partition $\la^{\prime}=(\la_{1}^{\prime}, \la_{2}^{\prime},\dots)$
so that $\la_{j}^{\prime}$ is the number of boxes in the $j$-th column of
the Young diagram $\la$. If $\alpha=(i, j)$
is a box of $\la,$ then the corresponding hook length
is defined as $h(\al)=\la_{i}+\la_{j}^{\prime}-i-j+1$ and
the content is $c(\al)=j-i$.
The {\em Schur element} $c_{\la}$ is given by
the formula
\beql{schel}
c_{\la}=\prod_{\al\in\la} q^{c(\al)}\tss [h(\al)]_q;
\eeq
here and below we use standard notation for $q$-integers,
\beql{qnum}
[k]_q = \frac{q^k -q^{-k}}{q-q^{-1}}\Fand [k]_q!=[1]_q\dots [k]_q.
\eeq

A $\la$-tableau with entries in
$\{1,\dots,n\}$ which are filled in the boxes of $\la$ bijectively,
is called {\em standard},
if its entries strictly increase along the rows and down the columns.

The Hecke algebra $\Hc_n$ is semisimple and isomorphic to the direct sum of matrix algebras
\beql{isomhecke}
\Hc_n \cong \bigoplus_{\la \vdash n} \Mat_{f_{\la}}(\CC).
\eeq
The diagonal matrix units $e_{\Uc}=e_{\Uc \Uc}\in \Mat_{f_{\la}}(\CC)$
parameterized by the standard tableaux of shape $\la$
are {\em primitive idempotents} of $\Hc_n$.
The {\em fusion procedure} originating in \cite{c:sb} and \cite{c:ni} provides an
explicit expression for $e_{\Uc}$ in terms of a rational function
in $n$ variables with values in $\Hc_n$; see \cite{n:mh} for detailed proofs.
A simple version of the procedure was developed in \cite{imo:ih}. Here we give
its brief review in the form similar the symmetric group version of \cite{m:fp}.
It is this approach which we will then
extend to the Hecke--Clifford superalgebra.

For all pairs $1\leqslant a<b\leqslant n$
introduce the corresponding rational functions $\vp_{ab}(u,v)$
in variables $u,v$ with values in $\Hc_n$
by setting
\beql{phiab}
\vp_{ab}(u,v)=1-\frac{\ve\tss u}{u-v}\ts
T_{b-1}^{-1}\dots T_a^{-1}\ts T_{a+1}\dots T_{b-1}.
\eeq
Now take $n$ complex variables $u_1,\dots,u_n$
and set
\beql{phiuuuh}
\Phi(u_1,\dots,u_n)=
\prod_{1\leqslant a<b\leqslant n}\vp^{}_{ab}(u_a,u_b),
\eeq
where the product is taken in the reverse lexicographic order on the pairs $(a,b)$:
\beql{ord}
(1,2)\prec (1,3)\prec (2,3)\prec \dots \prec (1,n)\prec\dots\prec (n-1,n).
\eeq

Given a standard
$\la$-tableau $\Uc$ and $a\in\{1,\dots,n\}$, we set
$\si_a=q^{2\tss(j-i)}$ if the element $a$ of $\Uc$ occupies
the box in row $i$ and column $j$.

\bth\label{thm:fphecke}
The consecutive evaluations are well-defined with the value given by
\beql{theor}
\Phi(u_1,\dots,u_n)
\big|_{u_1=\si_1} \big|_{u_2=\si_2} \dots \big|_{u_n=\si_n}=c_{\la'}\ts e_{\Uc},
\eeq
where $c_{\la'}$ is the Schur element associated with $\la'$.
\eth

\bpf
We argue by induction on $n$ and note that the relation holds trivially for $n=1$.
Hence
using
the induction hypothesis and setting $u=u_n$, for $n\geqslant 2$ we get
\ben
\Phi(u_1,\dots,u_n)\big|_{u_1=\si_1}\dots \big|_{u_{n-1}=\si_{n-1}}
=c_{\mu'}\ts e^{}_{\Vc}\ts
\vp^{}_{1n}(\si_1,u)\dots \vp^{}_{n-1,n}(\si_{n-1},u),
\een
where the standard tableau $\Vc$ is obtained from $\Uc$ by
removing the box occupied by $n$ and $\mu$ is the shape of $\Vc$.
Hence, it suffices to show that the following evaluation is well-defined with the value
given by
\beql{evah}
e^{}_{\Vc}\ts
\vp^{}_{1n}(\si_1,u)\dots \vp^{}_{n-1,n}(\si_{n-1},u)\big|_{u=\si_n}
=\frac{c_{\la'}}{c_{\mu'}}\ts e^{}_{\Uc}.
\eeq
To this end, introduce counterparts of the functions \eqref{phiab} by setting
\beql{phiba}
\vp_{ba}(u,v)=1-\frac{\ve\tss v}{u-v}\ts
T_{b-1}^{-1}\dots T_{a+1}^{-1}\ts T_{a}\dots T_{b-1}
\eeq
for $1\leqslant a<b\leqslant n$. Furthermore, recall that
the pairwise commuting {\em Jucys--Murphy elements} $x_1,\dots,x_n$ of $\Hc_{n}$
are defined by $x_1=1$ and
\beql{jmind}
x_{a+1} = T_a\dots T_2\tss T_1^2\tss T_2\dots T_a
\qquad\text{for}\quad a=1,\dots,n-1;
\eeq
see \cite{c:ni}, \cite{dj:bi}. They can be used to express the primitive idempotents
by an inductive rule. As with the case of the symmetric group
considered in \cite{m:fp}, the following rational function
in $u$
is regular at $u=\si_n$ with the value
\beql{idem1}
\frac{u - \si_n}{u - x_n}\ts e_{\Vc}\Big|_{u=\si_n}=e_{\Uc}.
\eeq
We have the properties
\beql{jmev}
x_a \ts e_{\Uc} =  e_{\Uc} \ts x_a = \si_a\ts e_{\Uc}
\qquad\text{for}\quad a=1,\dots,n.
\eeq

\ble\label{lem:prjmh}
We have the relation
\beql{redeh}
e^{}_{\Vc}\ts \vp^{}_{n,n-1}(u,\si_{n-1})\dots \vp^{}_{n1}(u,\si_{1})
=\frac{u-x_n}{u-1}\ts e^{}_{\Vc}.
\eeq
\ele

\bpf
It will be convenient to prove more general relations in $\Hc_n$,
\beql{redegenh}
e^{}_{\Vc}\ts \vp^{}_{n,r}(u,\si_{r})\dots \vp^{}_{n1}(u,\si_{1})
=\frac{u-x_{r n}}{u-1}\ts e^{}_{\Vc},
\eeq
where $r=1,\dots,n-1$ and $\Vc$ now denotes a standard tableau
with $r$ boxes, whereas
\ben
x_{rn}=T_{n-1}^{-1}\dots T_{r+1}^{-1}\ts T_{r}\dots T_2 T_{1}^2 T_2\dots T_r\ts T_{r+1}\dots T_{n-1}.
\een
Note that $x_{n}=x_{n-1,n}$.
Since $x_{rn}$ commutes with $x_1,\dots,x_r$, it also commutes with $e^{}_{\Vc}$.
It is clear that \eqref{redegenh} holds for $r=1$. To proceed by
the induction on $r$, suppose that $r\geqslant 2$.
Let $\Wc$ be the tableau obtained from $\Vc$ by
removing the box occupied by $r$. Write $e^{}_{\Vc}=e^{}_{\Vc}e^{}_{\Wc}$,
and observe that $e^{}_{\Wc}$ commutes with the factor $\vp^{}_{n,r}(u,\si_{r})$. Hence,
by the induction hypothesis, the left hand side of \eqref{redegenh}
can be written as
\beql{evpnr}
e^{}_{\Vc}\ts \vp^{}_{n,r}(u,\si_{r})\ts \frac{u-x_{r-1,n}}{u-1}\ts e^{}_{\Wc}
=e^{}_{\Vc}\ts \vp^{}_{n,r}(u,\si_{r})\ts \frac{u-x_{r-1,n}}{u-1}.
\eeq
Write $\vp^{}_{n,r}(u,\si_{r})$ in the form
\ben
\vp^{}_{n,r}(u,\si_{r})=\frac{u-\si_r\ts
T_{n-1}^{-1}\dots T_{r+1}^{-1}\ts T^2_{r}\ts T_{r+1}\dots T_{n-1}}{u-\si_r}.
\een
It remains to verify the identity
\ben
e^{}_{\Vc}\ts(u-\si_r\ts
T_{n-1}^{-1}\dots T_{r+1}^{-1}\ts T^2_{r}\ts T_{r+1}\dots T_{n-1})(u-x_{r-1,n})
=e^{}_{\Vc}\ts (u-\si_r)(u-x_{rn}).
\een
The constant terms of the quadratic polynomials in $u$ clearly agree.
Since $e^{}_{\Vc}\ts\si_r=e^{}_{\Vc}\ts x_r$, the required equality of the coefficients
of $-u$ reads
\ben
x_r\ts T_{n-1}^{-1}\dots T_{r+1}^{-1}\ts T^2_{r}\ts T_{r+1}\dots T_{n-1}+x_{r-1,n}=x_r+x_{rn}.
\een
However, this is equivalent to the relation
\ben
(x_{r}+x_{r+1})\ts T_r=T_r\ts (x_{r}+x_{r+1})
\een
which holds due to the properties of the Jucys--Murphy elements,
\ben
T_r\ts x_r=x_{r+1}\ts T_r-\ve\ts x_{r+1}\Fand T_r\ts x_{r+1}=x_{r}\ts T_r+\ve\ts x_{r+1},
\een
thus completing the proof of the lemma.
\epf

Returning to the proof of the theorem,
use the relations
\ben
\vp_{ab}(u,v)\ts \vp_{ba}(v,u)=A(u,v)
\een
which hold for $a\ne b$,
where $A(u,v)$ is defined by
\beql{auv}
A(u,v)=\frac{(u-q^2\tss v)(u-q^{-2}\tss v)}{(u-v)^2}.
\eeq
Hence, we derive from Lemma~\ref{lem:prjmh} that
\ben
e^{}_{\Vc}\ts
\vp^{}_{1n}(\si_1,u)\dots \vp^{}_{n-1,n}(\si_{n-1},u)
=\frac{u-1}{u-x_n}\ts e^{}_{\Vc}\ts A(u,\si_1)\dots A(u,\si_{n-1}),
\een
and the right hand side can be written as
\ben
\frac{u-\si_n}{u-x_n}\ts e^{}_{\Vc}\cdot
\frac{u-1}{u-\si_n}\ts A(u,\si_1)\dots A(u,\si_{n-1}).
\een
The required relation \eqref{evah} follows from the claim that both the following
evaluations are well-defined with the first factor evaluated by \eqref{idem1}, while
the formula
\beql{ratfnuh}
\frac{u-1}{u-\si_n}\ts A(u,\si_1)\dots A(u,\si_{n-1})\big|_{u=\si_n}=\frac{c_{\la'}}{c_{\mu'}}
\eeq
was already verified in \cite[Lemma~3.2]{imo:ih}.
\epf

To make connections with the versions of the fusion procedure of
\cite{c:ni}, \cite{imo:ih} and \cite{n:mh}, consider the rational functions
\ben
T_a(u,v)=T_a+\frac{\ve}{u^{-1}v-1},\qquad a=1,\dots,n-1,
\een
used therein. They satisfy the Yang--Baxter equation in the form
\beql{ybebr}
T_a(u,v) \, T_{a+1}(u,w) \, T_a(v,w) =
T_{a+1}(v,w) \, T_a(u,w) \, T_{a+1}(u,v).
\eeq

\bpr\label{prop:equivph}
The rational function \eqref{phiuuuh} can be written as
\beql{phil}
\Phi(u_1,\dots,u_n)=\prod_{1\leqslant a<b\leqslant n}
T_{b-a}(u_a, u_b)\cdot T_{\om^{}_n}^{-1},
\eeq
where the product is taken in the order \eqref{ord}, and
\beql{phir}
\Phi(u_1,\dots,u_n)=T_{\om^{}_n}^{-1}\cdot
\prod_{1\leqslant a<b\leqslant n}
T_{b-a}(u_a, u_b),
\eeq
where the product is taken in the order reverse to \eqref{ord}.
\epr

\bpf
It is clear that the formula \eqref{phil} is an immediate consequence of the relations
\ben
T_{\om^{}_{a-1}}\ts T_{a-1}(u_1,u_a)\dots T_{1}(u_{a-1},u_a)\ts T_{\om^{}_a}^{-1}
=\vp^{}_{1a}(u_1,u_a)\dots \vp^{}_{a-1,a}(u_{a-1},u_a)
\een
for $a=2,\dots,n$, which we will now verify. Indeed, the right hand side equals
\ben
T_{a-1}^{-1}\dots T_1^{-1}\ts T_{1}(u_1,u_a)\dots T_{a-1}(u_{a-1},u_a).
\een
Since $T_1\dots T_{a-1}\ts T_{\om^{}_{a-1}}=T_{\om^{}_a}$, the desired relations
follow from the properties
\beql{conle}
T_{\om^{}_a}\ts T_b(u,v)\ts T_{\om^{}_a}^{-1}=T_{a-b}(u,v)
\eeq
for $b=1,\dots,a-1$.

Furthermore, by a repeated application of \eqref{ybebr}, we can bring
the right hand side of \eqref{phil} to the form
\ben
\prod_{1\leqslant a<b\leqslant n}
T_{n-b+a}(u_a, u_b)\cdot T_{\om^{}_n}^{-1},
\een
where the product is taken in the order reverse to \eqref{ord}. Now \eqref{phir} follows from
\eqref{conle} with $a=n$.
\epf

By Proposition~\ref{prop:equivph}, the fusion
formula of Theorem~\ref{thm:fphecke} is equivalent to \cite[Thm~3.3]{imo:ih}.
Note that although the functions $\vp_{ab}(u,v)$
do not satisfy the Yang--Baxter equation,
their use makes the fusion formula \eqref{theor} analogous to
its symmetric group counterpart in \cite{m:fp}. We will introduce Hecke--Clifford analogues
of these functions in the next section to derive a version of the fusion
formula analogous to Theorem~\ref{thm:fphecke}.

\section{Fusion procedure for $\HC_n$}
\label{sec:fp}

As introduced by Olshanski~\cite{o:qu}, the
{\em Hecke--Clifford superalgebra} $\HC_n=\HC_n(q)$ is generated by even elements $T_1,\dots,T_{n-1}$
and odd elements $c_1,\dots,c_n$ subject to the defining relations \eqref{defhecke}
together with
\ben
c_a^2=-1,\qquad c_ac_b=-c_bc_a,\qquad a\ne b,
\een
and
\ben
\bal
T_ac_a&=c_{a+1}T_a,\qquad T_ac_{a+1}=c_aT_a-\ve(c_a-c_{a+1}),\\[0.4em]
T_ac_b&=c_{b}T_a,\qquad b\ne a, a+1.
\eal
\een
The subalgebra of $\HC_n$ generated by the elements $c_1,\dots,c_n$ is isomorphic
to the {\em Clifford superalgebra} $\Cl_n$.
As before, we assume that $q$ is a nonzero complex number which is not a root of unity.
This implies that the superalgebra $\HC_n$ is semisimple; see \cite{cw:rh}.
We follow \cite{jn:as} to define the
{\em Jucys--Murphy elements} in $\HC_n$ by $x_1=1$ and
\beql{jmse}
x_{a+1}=(T_a-\ve\ts c_ac_{a+1})\dots (T_1-\ve\ts c_1c_2)\ts T_1\dots T_a
\eeq
for $a=1,\dots,n-1$.

According to \cite{jn:as},
the irreducible representations of the superalgebra $\HC_n$ are
parameterized by strict partitions
$\la=(\la_1,\dots,\la_{\ell})$ of $n$ with $\la_1>\dots>\la_{\ell}>0$
and $\la_1+\dots+\la_{\ell}=n$. We will set $\ell(\la)=\ell$ to denote the length
of the partition and
write $\la\Vdash n$ to indicate that $\la$ is a strict partition of $n$.
A strict partition $\la$ is depicted by the shifted Young diagram obtained
from the usual Young diagram by shifting row $i$ to the right by $i-1$ unit boxes.

A ({\em shifted}) $\la$-{\em tableau} is obtained by writing the numbers $1,\dots,n$ bijectively
into the boxes of the shifted Young diagram $\la$.
Such a tableau
is called {\em standard} if its entries increase
from left to right in each row and from top to bottom in each column.
The entries in the boxes $(i,i)$ are called the {\em diagonal entries}.
As in \cite{kms:jm}, we will extend the set of tableaux by allowing any non-diagonal entry to
occur with a bar on it, as illustrated:
\ben
\begin{ytableau}
    1 & 2 & \ol{4} & \ol{5} & 8 & \ol{10}\\
    \none & 3 & \ol{6} & 9\\
    \none & \none & 7
\end{ytableau}
\ .
\een
We will use sans-serif letters $\asf, \bsf, \csf$, etc.
to denote arbitrary barred or unbarred
entries of $\Uc$. It will be convenient to use the sequence $\be_k$,
parameterized by integers $k$ which is
defined by
\beql{bek}
\be^{}_k=[k+1]_{q^2}-[k]_{q^2}+\ve\sqrt{[k+1]_{q^2}[k]_{q^2}},
\eeq
where we use notation \eqref{qnum}. Note that $\be_{-k}=\be_{k-1}$ and
\ben
\be^{-1}_k=[k+1]_{q^2}-[k]_{q^2}-\ve\sqrt{[k+1]_{q^2}[k]_{q^2}}.
\een
Given a standard barred tableau $\Uc$, introduce
the {\em $q$-contents} $\ka_{\asf}(\Uc)$
of any barred or unbarred entries $\asf$ of $\Uc$
by the formula
\beql{qcont}
\ka_{\asf}(\Uc)=\begin{cases}\be^{}_{\vs_{\ts a}}
\quad&\text{if $\asf=a$ is unbarred},\\[0.5em]
\be^{-1}_{\vs_{\ts a}}\quad&\text{if $\asf=\bar a$ is barred},
\end{cases}
\eeq
where $\vs_{\ts a}=j-i$ is the content of the box $(i,j)$ of $\la$
occupied by $a$ or $\bar a$ in $\Uc$.

For any standard barred tableau $\Uc$
introduce the element $e^{}_{\Uc}$ of $\HC_n$ by induction, setting
$e^{}_{\ts\tiny\young(1)}=1$ for the one-box tableau, and
\beql{murphyse}
e^{}_{\Uc}=e^{}_{\Vc}\ts \frac{(x_n-b_1)\dots
(x_n-b_p)}{(\ka-b_1)\dots (\ka-b_p)},
\eeq
where $\Vc$ is the barred tableau
obtained from $\Uc$ by removing the box $\al$ occupied by $n$ (resp. $\bar n$).
The shape of $\Vc$ is a shifted diagram which we denote by $\mu$ and
$b_1,\dots,b_p$ are the $q$-contents in all addable boxes of $\mu$ (barred and unbarred),
except for the entry $n$ (resp. $\bar n$), while $\ka$ is the $q$-content of
the entry $n$ (resp. $\bar n$).

The following proposition is a $q$-extension of the construction of \cite{kms:jm} for the
Sergeev superalgebra and it was already given in an equivalent form in the independent
work \cite{ls:sb}; cf. Remark~5.10 therein.

\bpr\label{prop:idemp}
All elements $e^{}_{\Uc}$ are idempotents in $\HC_n$. They are
pairwise orthogonal and form a decomposition of the identity:
\ben
e^{}_{\Uc}\tss e^{}_{\Vc}=\de^{}_{\Uc\Vc}\ts e^{}_{\Vc},\qquad
1=\sum_{\la\Vdash n}\ts\sum_{\sh(\Uc)=\la} e^{}_{\Uc}.
\een
Moreover,
\beql{xietsign}
x_a\ts e^{}_{\Uc}=e^{}_{\Uc}\ts x_a=\ka_{\asf}(\Uc)\ts e^{}_{\Uc},
\eeq
where $\ka_{\asf}(\Uc)$ is the $q$-content of the entry $\asf=a$ or $\asf=\bar a$  of\ \ts $\Uc$.
\epr

\bpf
As with the arguments of \cite{kms:jm} and \cite{ls:sb}, use
the explicit construction of $\HC_n$-modules as given in \cite[Sec.~6]{jn:as}.
For each strict partition $\la\Vdash n$
consider the module over $\HC_n$ afforded by the vector space
\beql{ula}
\wh U^{\la}=\bigoplus_{\sh(\Tc)=\la}\ts \Cl_n \tss v^{}_{\Tc}
\eeq
with the basis vectors $v^{}_{\Tc}$ associated with the (unbarred) standard $\la$-tableaux $\Tc$.
The Jucys--Murphy elements act on the basis vectors by
\beql{jmacti}
x_a\tss v^{}_{\Tc}=\ka_{a}(\Tc)\ts v^{}_{\Tc},
\eeq
for $a$ running over the set $\{1,\dots,n\}$. Since
the representation
\beql{uhat}
\wh U=\bigoplus_{\la\Vdash n}  \wh U^{\la}
\eeq
of the superalgebra $\HC_n$ is faithful, it suffices to verify
the desired properties of the elements of $\HC_n$ for their images
in the representation $\wh U$. This follows easily with the use
of the relations \eqref{jmacti} together with the properties
\beql{xaca}
x_ac_a=c_ax^{-1}_a,\qquad x_a\tss c_b=c_b\tss x_a,\qquad a\ne b,
\eeq
recorded in \cite[Sec.~3]{jn:as}.
\epf

Introduce rational functions in variables $u,v$ with values in $\HC_n$
by setting for $a<b$,
\beql{phiabhc}
\vp_{ab}(u,v)=T_{b-1}^{-1}\dots T_a^{-1}\ts \psi_{a}(u,v)\ts T_{a+1}\dots T_{b-1}
\eeq
and
\ben
\vp_{ba}(u,v)=T_{b-1}^{-1}\dots T_{a+1}^{-1}\ts \psi_{a}(u,v)\ts T_a\dots T_{b-1},
\een
where $\psi_a(u,v)$ is given by
\beql{psiuv}
\psi_a(u,v)=T_a+\frac{\ve}{u^{-1}v-1}+\frac{\ve\ts c_ac_{a+1}}{uv-1},
\eeq
as defined in \cite[Sec.~3]{jn:as}.
Take $n$ complex variables $u_1,\dots,u_n$
and set
\beql{phiuuuhc}
\Phi(u_1,\dots,u_n)=
\prod_{1\leqslant a<b\leqslant n}\vp^{}_{ab}(u_a,u_b),
\eeq
where the product is taken in the order \eqref{ord}.

Suppose that $\la\Vdash n$
and fix a standard barred $\la$-tableau $\Uc$.
For every $a=1,\dots,n$ set
$\ka_a=\ka_{\asf}(\Uc)$ if $\asf=a$ or $\asf=\bar a$ is the entry of $\Uc$
and consider the product
\ben
\ka^{}_{\Uc}:=\prod_{a=1}^n(1+\ka_a^{-1}).
\een
Furthermore, introduce a strict partition counterpart of the element \eqref{schel}
by setting
\beql{schelstr}
c_{\la}=\prod_{\al\in\la} \frac{1}{q^{c(\al)}+q^{-c(\al)}}\tss \prod_{a=1}^{\ell}
[\la_a]_q!\ts \prod_{1\leqslant i<j\leqslant \ell}\ts
\frac{[\la_i+\la_j]_q}{[\la_i-\la_j]_q},
\eeq
where $c(\al)=j-i$ is the content of the box $\al=(i,j)$.

The following is our main result providing an analogue of Theorem~\ref{thm:fphecke} for the Hecke--Clifford superalgebra
$\HC_n$. A dual version of this result will be proved in Theorem~\ref{thm:fpheckeclianti} below.

\bth\label{thm:fpheckecli}
The consecutive evaluations of the function \eqref{phiuuuhc}
are well-defined with the value given by
\beql{theorhc}
\Phi(u_1,\dots,u_n)
\big|_{u_1=\ka_1} \big|_{u_2=\ka_2} \dots \big|_{u_n=\ka_n}=c_{\la}\ts \ka^{}_{\Uc}\ts e^{}_{\Uc}.
\eeq
\eth

\bpf
Arguing as in the proof of Theorem~\ref{thm:fphecke}, use
induction on $n$. Setting $u=u_n$, for $n\geqslant 2$ we get
\ben
\Phi(u_1,\dots,u_n)\big|_{u_1=\ka_1}\dots \big|_{u_{n-1}=\ka_{n-1}}
=c_{\mu}\ts \ka^{}_{\Vc}\ts e^{}_{\Vc}\ts
\vp^{}_{1n}(\ka_1,u)\dots \vp^{}_{n-1,n}(\ka_{n-1},u),
\een
where the standard tableau $\Vc$ is obtained from $\Uc$ by
removing the box occupied by $n$ or $\bar n$ and $\mu$ is the shape of $\Vc$.
Hence, it suffices to show that the following evaluation is well-defined with the value
given by
\beql{evahhc}
e^{}_{\Vc}\ts
\vp^{}_{1n}(\ka_1,u)\dots \vp^{}_{n-1,n}(\ka_{n-1},u)\big|_{u=\ka_n}
=\frac{c_{\la}\ts \ka^{}_{\Uc}}{c_{\mu}\ts \ka^{}_{\Vc}}\ts e^{}_{\Uc}.
\eeq

\ble\label{lem:prjm}
We have the relation
\beql{redehc}
e^{}_{\Vc}\ts \vp^{}_{n,n-1}(u,\ka_{n-1})\dots \vp^{}_{n1}(u,\ka_{1})
=\frac{u-x_n}{u-1}\ts e^{}_{\Vc},
\eeq
where the barred tableau $\Vc$ is obtained from $\Uc$ by
removing the box occupied by $n$ or $\bar n$.
\ele

\bpf
We will modify the argument used to prove Lemma~\ref{lem:prjmh} to derive more general relations in $\HC_n$,
\beql{redegen}
e^{}_{\Vc}\ts \vp^{}_{n,r}(u,\ka_{r})\dots \vp^{}_{n1}(u,\ka_{1})
=\frac{u-x_{rn}}{u-1}\ts e^{}_{\Vc},
\eeq
where $r=1,\dots,n-1$ and $\Vc$ now denotes a standard barred tableau
with $r$ boxes. The numbers $\ka_a=\ka_{\asf}(\Vc)$ are the $q$-contents of the entries of $\Vc$,
whereas
\ben
x_{rn}=T_{n-1}^{-1}\dots T_{r+1}^{-1}\ts (T_r-\ve\ts c_rc_{r+1})\dots (T_1-\ve\ts c_1c_2)\ts T_{1}\dots T_r\ts T_{r+1}\dots T_{n-1},
\een
so that
\beql{xrn}
x_{rn}=T_{n-1}^{-1}\dots T_{r+1}^{-1}\ts x_{r+1}\ts T_{r+1}\dots T_{n-1}.
\eeq
Since $x_{rn}$ commutes with $x_1,\dots,x_r$, it also commutes with $e^{}_{\Vc}$.
It is easily seen that \eqref{redegen} holds for $r=1$. We proceed by
induction on $r$, and note that for $r\geqslant 2$ the left hand side of \eqref{redegen}
can be written as
\beql{evpnrgen}
e^{}_{\Vc}\ts \vp^{}_{n,r}(u,\ka_{r})\ts \frac{u-x_{r-1,n}}{u-1}.
\eeq
Write $\vp^{}_{n,r}(u,\ka_{r})$ in the form
\ben
\vp^{}_{n,r}(u,\ka_{r})=T_{n-1}^{-1}\dots T_{r+1}^{-1}
\Bigg(1-\Big(\frac{\ve\ts\ka_r}{u-\ka_r}-\frac{\ve\ts c_rc_{r+1}}{u\tss\ka_r-1}\Big)\tss T_r\Bigg)
T_{r+1}\dots T_{n-1}.
\een
Taking \eqref{xrn} into account, we come
to verify the identity
\beql{keyid}
e^{}_{\Vc}\ts \Bigg(1-\Big(\frac{\ve\ts\ka_r}{u-\ka_r}-\frac{\ve\ts c_rc_{r+1}}{u\tss\ka_r-1}\Big)\tss T_r\Bigg)
\big(u-T_{r}^{-1} x_r\ts T_r\big)
=e^{}_{\Vc}\ts (u-x_{r+1}).
\eeq
Since $e^{}_{\Vc}\ts\ka_r=e^{}_{\Vc}\ts x_r$, using \eqref{xaca}
we find that
\ben
e^{}_{\Vc}\ts c_rc_{r+1}\tss x_r=e^{}_{\Vc}\ts c_rc_{r+1}\tss \ka^{-1}_r.
\een
Hence, we can write the expression
\ben
e^{}_{\Vc}\ts\Big(\frac{\ve\ts\ka_r}{u-\ka_r}-\frac{\ve\ts c_rc_{r+1}}{u\tss\ka_r-1}\Big)\tss T_r
\big(u-T_{r}^{-1} x_r\ts T_r\big)
\een
in the form
\ben
\frac{\ve}{u-\ka_r}\ts e^{}_{\Vc}\ts (u\ts x_r\ts T_r-\ka_r\ts x_r\ts T_r)
-\frac{\ve}{u\tss\ka_r-1}\ts e^{}_{\Vc}\ts c_rc_{r+1}\ts (u\ts T_r-\ka^{-1}_r\ts T_r)
\een
which simplifies to
\ben
\ve\ts e^{}_{\Vc}\ts x_r\ts T_r-\ve\ts e^{}_{\Vc}\ts c_rc_{r+1}\tss x_r\ts T_r.
\een
Recalling that $T_r^{-1}=T_r-\ve$, we derive \eqref{keyid}, thus
completing the proof of the lemma.
\epf

Returning to the proof of the theorem,
use the relations
\ben
\vp_{ab}(u,v)\ts \vp_{ba}(v,u)=A(u,v)
\een
which hold for $a\ne b$, as implied by \cite[Lemma~4.1]{jn:as}, where
\ben
A(u,v)=1-\ve^2\ts \Bigg(\frac{u^{-1}v}{(u^{-1}v-1)^2}+\frac{u\tss v}{(u\tss v-1)^2}\Bigg).
\een
This expression can also be written as
\beql{auvnew}
A(u,v)=\frac{(u+u^{-1})^2+(v+v^{-1})^2-(\ve^2+2)(u+u^{-1})(v+v^{-1})+4\tss \ve^2}
{(u+u^{-1}-v-v^{-1})^2}.
\eeq
We derive from Lemma~\ref{lem:prjm} that
\ben
e^{}_{\Vc}\ts
\vp^{}_{1n}(\ka_1,u)\dots \vp^{}_{n-1,n}(\ka_{n-1},u)
=\frac{u-1}{u-x_n}\ts e^{}_{\Vc}\ts A(u,\ka_1)\dots A(u,\ka_{n-1}),
\een
and the right hand side can be written as
\ben
\frac{u-\ka_n}{u-x_n}\ts e^{}_{\Vc}\cdot
\frac{u-1}{u-\ka_n}\ts A(u,\ka_1)\dots A(u,\ka_{n-1}).
\een
The required relation  will follow from the claim that both the following
evaluations are well-defined with the values given by
\beql{ratfu}
\frac{u-\ka_n}{u-x_n}\ts e^{}_{\Vc}\big|_{u=\ka_n}=e^{}_{\Uc}
\eeq
and
\beql{ratfnu}
\frac{u-1}{u-\ka_n}\ts A(u,\ka_1)\dots A(u,\ka_{n-1})\big|_{u=\ka_n}
=\frac{c_{\la}\ts \ka^{}_{\Uc}}{c_{\mu}\ts \ka^{}_{\Vc}}.
\eeq
Relation \eqref{ratfu} is a Hecke--Clifford counterpart of
\eqref{idem1} and verified by the same argument as for the symmetric group in \cite{m:fp};
cf. the proof of \cite[Theorem~4.1]{kms:jm}.
It relies on the property of the idempotents
\ben
e^{}_{\Vc}=\sum_{\Vc\nearrow \Wc} e^{}_{\Wc},
\een
where $\Vc\sclnearrow \Wc$ means that
the barred tableau $\Wc$ is obtained from $\Vc$
by adding one box with the entry $n$ or $\bar n$.

To verify \eqref{ratfnu}, note that the rational function $A(u,v)$
is stable under the replacement $v\mapsto v^{-1}$ and so we may
take $\Vc$ to be unbarred. Moreover,
the product
depends only on the shape $\mu$ of $\Vc$ so we may choose a particular
tableau $\Vc$ for its evaluation. We will take
the row tableau
obtained by filling in the boxes of $\mu$ with
the numbers $1,\dots,n-1$ by consecutive rows from left to right
in each row.

Observe that the relation \eqref{auvnew} implies
a factorized expression
\ben
A(u,\be_{p})=\frac{(u+u^{-1}-\be^{}_{p-1}-\be^{-1}_{p-1})\tss(u+u^{-1}-\be^{}_{p+1}-\be^{-1}_{p+1})}
{(u+u^{-1}-\be^{}_{p}-\be^{-1}_{p})^2}.
\een
Noting that $\be_{-1}=\be_0=1$, for $p\geqslant 1$ we get
the formula
\ben
A(u,\be_0)\dots A(u,\be_{p-1})=\frac{(u-\be^{}_{p})(u-\be^{-1}_{p})}{(u-\be^{}_{p-1})(u-\be^{-1}_{p-1})}.
\een
Therefore, introducing the parts of $\mu$ by $\mu=(\mu_1,\dots,\mu_m)$ with $m=\ell(\mu)$, we derive
\ben
A(u,\ka_1)\dots A(u,\ka_{n-1})=\prod_{k=1}^m\ts \frac{(u-\be^{}_{\mu_k})(u-\be^{-1}_{\mu_k})}{(u-\be^{}_{\mu_k-1})(u-\be^{-1}_{\mu_k-1})}.
\een

Suppose first that
$m=\ell(\la)$ so that the tableau $\Uc$ is obtained from $\Vc$ by adding a
box with the entry $n$ or $\bar n$ at the row $r\in\{1,\dots,m\}$. Then
$\ka_n=\be^{}_{\mu_r}$ or $\ka_n=\be^{-1}_{\mu_r}$, respectively,
while $\mu_{r-1}\geqslant \mu_r+2$. Therefore, the left hand side of \eqref{ratfnu}
is well-defined with the respective values
\ben
\frac{(\be^{}_{\mu_r}-1)(\be^{}_{\mu_r}-\be^{-1}_{\mu_r})}
{(\be^{}_{\mu_r}-\be^{}_{\mu_r-1})(\be^{}_{\mu_r}-\be^{-1}_{\mu_r-1})}
\prod_{k=1, \ k\ne r}^m\ts \frac{(\be^{}_{\mu_r}-\be^{}_{\mu_k})(\be^{}_{\mu_r}-\be^{-1}_{\mu_k})}
{(\be^{}_{\mu_r}-\be^{}_{\mu_k-1})(\be^{}_{\mu_r}-\be^{-1}_{\mu_k-1})}
\een
and
\ben
\frac{(\be^{-1}_{\mu_r}-1)(\be^{-1}_{\mu_r}-\be^{}_{\mu_r})}
{(\be^{-1}_{\mu_r}-\be^{}_{\mu_r-1})(\be^{-1}_{\mu_r}-\be^{-1}_{\mu_r-1})}
\prod_{k=1, \ k\ne r}^m\ts \frac{(\be^{}_{\mu_r}-\be^{}_{\mu_k})(\be^{}_{\mu_r}-\be^{-1}_{\mu_k})}
{(\be^{}_{\mu_r}-\be^{}_{\mu_k-1})(\be^{}_{\mu_r}-\be^{-1}_{\mu_k-1})}.
\een
In both cases, the expressions can be written in a uniform way as
\ben
\big(\ka_n^{-1}+1\big)\ts\frac{[\mu_r+1]_q}{q^{\mu_r}+q^{-\mu_r}}\ts\prod_{k=1,\ k\ne r}^m
\frac{[\mu_r-\mu_k]_q\ts [\mu_r+\mu_k+1]_q}{[\mu_r+\mu_k]_q\ts [\mu_r-\mu_k+1]_q}.
\een
It is easily seen that this
coincides with the right hand side of
\eqref{ratfnu}, thus completing proof in the case under consideration.

In the remaining case with $m=\ell(\la)-1$, the tableau $\Uc$
is obtained from $\Vc$ by adding the
box $(\ell(\la),\ell(\la))$ with the entry $n$. This is possible when $\mu_m\geqslant 2$.
In this case, $\ka_n=1$ and
the left hand side of \eqref{ratfnu}
is well-defined with the value
\ben
\prod_{k=1}^m\ts \frac{[\mu_k+1]_q}{[\mu_k-1]_q},
\een
which coincides with the right hand side of \eqref{ratfnu}.
\epf

\bre\label{rem:qone}
In the limit $q\to 1$, the value of the constant factor $c_{\la}\ts \ka^{}_{\Uc}$ in \eqref{theorhc}
equals $n!/g_{\la}$, where $g_{\la}$ is the number of standard $\la$-tableaux
found by the Schur formula
\ben
g_{\la}=\frac{n!}{\la_1!\dots \la_{\ell}!}\ts \prod_{1\leqslant i<j\leqslant \ell}\ts
\frac{\la_i-\la_j}{\la_i+\la_j}.
\een
This limit value occurs as the constant factor in the fusion formula for the Sergeev superalgebra
in \cite[Thm~4.1]{kms:jm}. Moreover, by setting $u=q^{2x}$ and $v=q^{2y}$ for variables $x,y$,
we find that the limit value of the function \eqref{phiabhc} in the Sergeev superalgebra
is given by
\ben
\vp_{ab}(u,v)\to 1-\frac{(a,b)}{x-y}+\frac{(a,b)\ts c_ac_b}{x+y},\qquad \text{as}\quad q\to 1.
\een
For the $q$-contents defined in \eqref{qcont} we have the respective limits
\ben
\frac{\ka_{\asf}(\Uc)-1}{\ve}\to \pm\ts\sqrt{\vs_{\ts a}(\vs_{\ts a}+1)}.
\een
Hence, applying the limit $q\to 1$ to the fusion formula of Theorem~\ref{thm:fpheckecli},
we derive the fusion formula for the Sergeev superalgebra
in \cite[Thm~4.1]{kms:jm}.
\qed
\ere

We now give an equivalent formula for the rational function
\eqref{phiuuuhc} in terms of the functions \eqref{psiuv}, similar to
Proposition~\ref{prop:equivph}.

\bpr\label{prop:equivphc}
The rational function \eqref{phiuuuhc} can be written as
\beql{phirhc}
\Phi(u_1,\dots,u_n)=T_{\om^{}_n}^{-1}\cdot
\prod_{1\leqslant a<b\leqslant n}
\psi_{b-a}(u_a, u_b)
\eeq
where the product is taken in the order on the pairs $(a,b)$ reverse to \eqref{ord}.
\epr

\bpf
We start by proving by induction the formula
\beql{phirhcin}
\Phi(u_1,\dots,u_n)=T_{\om^{}_n}^{-1}\cdot
\prod_{1\leqslant a<b\leqslant n}
\psi_{n+a-b}(u_a, u_b),
\eeq
where the product is taken in the order \eqref{ord} on the pairs $(a,b)$.
The formula holds trivially for $n=1$ so suppose that $n\geqslant 2$.
By the induction hypothesis, we can write
\ben
\Phi(u_1,\dots,u_n)=T_{\om^{}_{n-1}}^{-1}\cdot
\prod_{1\leqslant a<b\leqslant n-1}
\psi_{n+a-b-1}(u_a, u_b)\cdot\vp_{1,n}(u_1,u_n)\dots \vp_{n-1,n}(u_{n-1},u_n).
\een
On the other hand, the definition \eqref{phiabhc} implies
\ben
\vp_{1,n}(u_1,u_n)\dots \vp_{n-1,n}(u_{n-1},u_n)=T_{n-1}^{-1}\dots T_1^{-1}\ts
\psi_{1}(u_1,u_n)\dots \psi_{n-1}(u_{n-1},u_n).
\een
The formula \eqref{phirhcin} now follows from the easily verified relations
$T_{\om^{}_n}=T_1\dots T_{n-1}\ts T_{\om^{}_{n-1}}$ and
\ben
\psi_k(u,v)\ts T_{n-1}^{-1}\dots T_1^{-1}=T_{n-1}^{-1}\dots T_1^{-1}\ts \psi_{k+1}(u,v)
\een
for $k=1,\dots,n-2$.
The expression \eqref{phirhc} is equivalent to \eqref{phirhcin}
as already pointed out in \cite[Proposition~5.5]{jn:as}
due to the relations satisfied
by the functions \eqref{psiuv}; see \cite[Lemma~4.1]{jn:as}:
\ben
\psi_{k}(u,v)\ts \psi_{l}(w,y)=\psi_{l}(w,y)\ts \psi_{k}(u,v),\qquad |k-l|\geqslant 2,
\een
and
\ben
\psi_{k}(u,v)\ts \psi_{k+1}(w,v)\ts \psi_{k}(w,u)=\psi_{k+1}(w,u)\ts \psi_{k}(w,v)\ts\psi_{k+1}(u,v);
\een
the latter is known as the {\em Yang--Baxter equation}.
\epf

Introduce elements of the Hecke--Clifford superalgebra $\HC_n$ by setting $\BT_a=T_a-\ve\ts c_a\tss c_{a+1}$.
As was pointed out in \cite[Sec.~2-i]{bk:hc}, the mapping
\beql{aniinv}
T_a\mapsto \BT_a,\qquad c^{}_b\mapsto c^{}_b
\eeq
for all $a=1,\dots,n-1$ and $b=1,\dots,n$, defines an anti-involution of $\HC_n$.
Denote by $\oPhi\tss(u_1,\dots,u_n)$ the image of the rational function \eqref{phiuuuhc}
under this anti-involution. Explicitly,
\beql{phiuuuhcanti}
\oPhi\tss(u_1,\dots,u_n)=
\prod_{1\leqslant a<b\leqslant n}\ovp^{}_{ab}(u_a,u_b),
\eeq
where the product is taken in the order reverse to \eqref{ord}, while
\ben
\ovp^{}_{ab}(u,v)=\BT_{b-1}\dots \BT_{a+1}\ts \opsi_{a}(u,v)\ts \BT_{a}^{-1}\dots \BT_{b-1}^{-1}
\een
and
\ben
\opsi_{a}(u,v)=T_a+\frac{\ve}{u^{-1}v-1}-\frac{\ve\ts c_ac_{a+1}\ts u\tss v}{u\tss v-1}.
\een

The following is a dual version of Theorem~\ref{thm:fpheckecli} which holds under the same assumptions.

\bth\label{thm:fpheckeclianti}
The consecutive evaluations of the function \eqref{phiuuuhcanti}
are well-defined with the value given by
\beql{theorhcanti}
\oPhi(u_1,\dots,u_n)
\big|_{u_1=\ka_1} \big|_{u_2=\ka_2} \dots \big|_{u_n=\ka_n}=c_{\la}\ts \ka^{}_{\Uc}\ts e^{}_{\Uc}.
\eeq
\eth

\bpf
Observe that the Jucys--Murphy elements \eqref{jmse} are stable under the
anti-involution \eqref{aniinv} and so are the idempotents defined in
\eqref{murphyse}.
Using the same induction argument as in the proof of Theorem~\ref{thm:fpheckecli},
we come to verifying that the following evaluation is well-defined with the value
given by
\ben
\ovp^{}_{n-1,n}(\ka_{n-1},u)\dots
\ovp^{}_{1n}(\ka_1,u)\ts e^{}_{\Vc}\big|_{u=\ka_n}
=\frac{c_{\la}\ts \ka^{}_{\Uc}}{c_{\mu}\ts \ka^{}_{\Vc}}\ts e^{}_{\Uc}.
\een
This formula is immediate from \eqref{evahhc} by the application of
the
anti-involution \eqref{aniinv}.
\epf

We also have a dual version of Proposition~\ref{prop:equivphc} which is immediate
by the application of the anti-involution \eqref{aniinv}.

\bpr\label{prop:equivphcanti}
The rational function \eqref{phiuuuhcanti} can be written as
\beql{phirhcanti}
\oPhi(u_1,\dots,u_n)=
\prod_{1\leqslant a<b\leqslant n}
\opsi_{b-a}(u_a, u_b)\cdot \BT_{\om^{}_n}^{-1},
\eeq
where the product is taken in the order \eqref{ord} on the pairs $(a,b)$,
and $\BT_{\om^{}_n}$ is the image of $T_{\om^{}_n}$ under the anti-involution \eqref{aniinv}.
\qed
\epr

As a concluding remark, we note that by Proposition~\ref{prop:equivph},
the rational function \eqref{phiuuuh} is stable
under the anti-involution of the Hecke algebra $\Hc_n$ identical on the generators $T_a$.
On the other hand, the rational functions \eqref{phiuuuhc} and \eqref{phiuuuhcanti}
are different already for $n= 2$. Nevertheless, Theorems~\ref{thm:fpheckecli}
and \ref{thm:fpheckeclianti} show that their evaluations at the $q$-contents coincide.

\section*{Declarations}

\subsection*{Competing interests}
The authors have no competing interests to declare that are relevant to the content of this article.

\subsection*{Acknowledgements}
Our work was supported by the Simons Foundation (grant MP-TSM-00002518 -- Jing),
the National Natural Science Foundation of China (grants 12471026 -- Liu and 12571026 -- Zhang),
and the Australian Research Council (grant DP240101572 -- Molev).

\subsection*{Availability of data and materials}
No data was used for the research described in the article.


\bigskip
\bigskip

\small

\noindent
N.J.:\newline
Department of Mathematics\\
North Carolina State University, Raleigh, NC 27695, USA\\
jing@ncsu.edu

\vspace{5 mm}

\noindent
I.K.:\newline
Shenzhen International Center for Mathematics\\
Southern University of Science and Technology, China\\
{\tt kashuba@sustech.edu.cn}

\vspace{5 mm}

\noindent
M.L.:\newline
School of Artificial Intelligence\newline
Jianghan University,
Wuhan 430056, Hubei, China\newline
{\tt ming.l1984@gmail.com}

\vspace{5 mm}

\noindent
A.M.:\newline
School of Mathematics and Statistics\newline
University of Sydney,
NSW 2006, Australia\newline
alexander.molev@sydney.edu.au

\vspace{5 mm}

\noindent
J.Z.:\newline
School of Mathematics and Statistics and Hubei Key Lab - Math. Sci.\newline
Central China Normal University, Wuhan 430079, Hubei, China\newline
{\tt jzhang@ccnu.edu.cn}

\end{document}